\documentclass[12pt]{amsart}

\usepackage{comment}
\usepackage[normalem]{ulem}
\usepackage[colorlinks,linkcolor=red,anchorcolor=blue,citecolor=blue]{hyperref} 
\usepackage{amscd} 
\usepackage{amsmath}
\usepackage{mathtools}
\usepackage{amssymb} 
\usepackage{amsthm}
\usepackage{dsfont}
\usepackage{bigdelim}
\usepackage{color} 
\usepackage{enumerate}
\usepackage{graphicx}
\usepackage{mathrsfs}
\usepackage{multirow}
\usepackage[all]{xy} 
\usepackage{color} 
\usepackage{enumitem}
\usepackage[dvipsnames]{xcolor}

\newtheorem{Theorem}{Theorem}[section] 
\newtheorem*{MetaTheorem}{Meta Theorem} 
\newtheorem{Corollary}[Theorem]{Corollary}
\newtheorem{Proposition}[Theorem]{Proposition}

\newtheorem{Lemma}[Theorem]{Lemma}
\theoremstyle{definition} 
\newtheorem{Definition}[Theorem]{Definition}
\newtheorem{Example}[Theorem]{Example}

\theoremstyle{remark}
\newtheorem{Remark}[Theorem]{Remark}

\newtheorem{Construction}[Theorem]{Construction}

\numberwithin{equation}{section}

\DeclareMathOperator{\Aut}{Aut}

\DeclareMathOperator{\GL}{GL}

\DeclareMathOperator{\Alb}{Alb}

\DeclareMathOperator{\cc}{c}
\DeclareMathOperator{\ch}{ch}

\DeclareMathOperator{\rk}{rk}

\newcommand{\R}{\mathbb{R}}
\newcommand{\Z}{\mathbb{Z}}

\newcommand{\C}{\mathbb{C}}
\newcommand{\Q}{\mathbb{Q}}

\newcommand{\mP}{\mathbb{P}}

\newcommand{\E}{\mathbb{E}}

\makeatletter
\newcommand*\bigcdot{\mathpalette\bigcdot@{.5}}
\newcommand*\bigcdot@[2]{\mathbin{\vcenter{\hbox{\scalebox{#2}{$\m@th#1\bullet$}}}}}
\makeatother

\title[Locally Constant Fibrations \& Positivity of Curvature]
{Locally Constant Fibrations \& Positivity of Curvature}

\author{Niklas M\"uller}
\address{Department of Mathematics, Universit\"at Duisburg-Essen,
Thea-Leymann-Str. 9, 45127 Essen, Germany.}
\email{{\tt niklas.mueller@uni-duisburg-essen.de}}

\date{\today}



\usepackage[capitalize,nameinlink]{cleveref}
\crefname{Lemma}{Lemma}{Lemmata}
\crefname{Definition}{Definition}{Definitions}
\crefname{Example}{Example}{Examples}
\crefname{Theorem}{Theorem}{Theorems}
\crefname{Corollary}{Corollary}{Corollaries}
\crefname{Proposition}{Proposition}{Propositions}
\crefname{Equation}{Eq.}{Eq.}
\crefname{IEEEeqnarray}{Eq.}{Eq.}
\crefname{Remark}{Remark}{Remarks}
\crefname{Reminder}{Reminder}{Reminders}
\crefname{figure}{Figure}{Figures}
\crefname{chapter}{Chapter}{Chapters}
\crefname{Construction}{Construction}{Constructions}
\crefname{align}{Eq.}{Eq.}
\crefname{align*}{Eq.}{Eq.}
\crefname{Conjecture}{Conjecture}{Conjectures}
\crefname{Question}{Question}{Questions}

\begin{document}
	\begin{abstract}
 
	\noindent
		Up to finite étale cover, any smooth complex projective variety $X$ with nef anti-canonical bundle is a holomorphic fibre bundle over a smooth projective variety with trivial canonical class (K-trivial variety for short) with locally constant transition functions. We show that this result is optimal by proving that any projective fibre bundle with locally constant transition functions over a $K$-trivial variety has a nef anti-canonical bundle. Moreover, we complement some results on the structure theory of varieties whose tangent bundle admits a singular Hermitian metric of positive curvature.
	\end{abstract}
    \maketitle
	\tableofcontents

	\section{Introduction}
	Starting with Mori's famous resolution of Hartshorne's conjecture, stating that projective space is the only complex manifold with ample tangent bundle, much work has been dedicated to the structure theory of varieties of positive or, more generally, non-negative curvature. The last years have seen many results of the form below published; some prime examples are \cite{DPS_ManifoldsWithNefTangentBundle, campanaDemaillyPeternell_SemiPositiveRicci, cao_NefAnticanonicalBundleII, cao_NefAnticanonicalBundleIII, DreulBianco_NumFlatLogTangentBundle, matsumuraWang_NefAnticanonicalBundle, hosonoIwaiMatsumura_PsefTangentBundle, Ejiri_varietiespositivecharacteristicnumerically, LOWYZ_FibrationsSTrictlyNefRelativeAC}:
	
	\begin{MetaTheorem}
		
		\noindent
		Let $X$ be a smooth projective variety whose tangent bundle/ anti-canonical bundle satisfies a suitable positivity condition $(\mathcal{P})$. Then there exists a finite étale cover $X'\rightarrow X$ and a smooth fibration 
        $$f\colon X'\rightarrow Y$$
        onto a smooth variety $Y$ with numerically trivial canonical class $\cc_1(K_Y) = 0$.
        Moreover, the general fibre $F$ of $f$ is rationally connected and satisifies $(\mathcal{P})$.
	\end{MetaTheorem}
	In fact, often much more can be said about $f$: typically it is a fibre bundle and often the transition functions can even be chosen to be constant. In the latter case we call $f$ a \emph{locally constant fibration}, see Definition \ref{ex:LocallyConstantFibration} for the precise definition.
	
	The first goal of this paper is to study the converse question: Given a locally constant fibration $f\colon X \rightarrow Y$ over a $K$-trivial base $Y$ such that the general fibre of $f$ satisfies $(\mathcal{P})$, is it true that $X$ satisfies $(\mathcal{P})$? In other words, we want to know whether the results obtained in \emph{loc.\ cit.\ }are sharp. Our main result is that this is true when $(\mathcal{P})$ is the nefness of the anti-canonical bundle:
	\begin{Theorem}
		Let $Y$ be a smooth projective variety with $\cc_1(Y) = 0$, let $F$ be a smooth, projective, rationally connected variety whose anti-canonical divisor $-K_F$ is nef and let $\rho\colon \pi_1(Y) \rightarrow \Aut(F)$ be any group homomorphism. Write
		\begin{align*}
		    X := (\widetilde{Y}\times F )/\pi_1(Y),
		\end{align*}
		where we denote by $\widetilde{Y}\rightarrow Y$ the universal cover of $Y$ and where the action of $\pi_1(Y)$ on $\widetilde{Y}\times F$ is given by $\gamma\colon \widetilde{Y}\times F \rightarrow \widetilde{Y}\times F,$ $(y, x)\mapsto (\gamma\cdot y, \rho(\gamma)(x))$ for any $\gamma\in \pi_1(Y)$. Then
		\begin{itemize}
		    \item[(1)] the smooth, compact, complex analytic variety $X$ is projective if and only if the image of $\pi_1(Y) \overset{\rho}{\rightarrow} \Aut(F)/\Aut^0(F)$ is finite.
			\item[(2)] In case $X$ is projective the anti-canonical divisor $-K_X$ is nef.
		\end{itemize}
		\label{IntroCharacterisationMfdsNefACBundle}
	\end{Theorem}
    Conversely, it is known by work of \cite{cao_NefAnticanonicalBundleIII} that, up to taking finite étale quotients, any smooth projective variety $X$ with a nef anti-canonical divisor arises in this way. We would like to emphasise that Theorem \ref{IntroCharacterisationMfdsNefACBundle} allows to easily construct many new non-trivial examples of varieties whose anti-canonical bundle is nef.

    The conclusion in Theorem \ref{IntroCharacterisationMfdsNefACBundle} is somewhat surprising to the author as the image of $\rho$ need not be compact; in particular, the action of $\pi_1(Y)$ on $F$ need not preserve any metric. Analogous statements hold for varieties with nef tangent bundle (Proposition \ref{characterisationMfdNefTangent}) and for mildly singular spaces (Theorem \ref{charMfdsNefALogCanBundle}). In case one asks for more positivity, the sharpness of the structural results obtained in the literature before is mostly immediate (Proposition \ref{charOtherPosNotions}). In the nef case however, to the best of the author's knowledge, apart from the partial results in \cite[Theorem 1.4]{matsumura_asymptoticBaseLoci}, these questions have not been previously addressed in the literature. 

    Note that the conclusion in Theorem \ref{IntroCharacterisationMfdsNefACBundle} is non-trivial and depends crucially on the assumption $c_1(Y) = 0$, c.f. Example \ref{ex:PositivityCanFail}. The main ingredient which makes Theorem \ref{IntroCharacterisationMfdsNefACBundle} work is the following result; in the smooth case it follows immediately from an observation due to Biswas \cite[Remark 3.7.(ii)]{biswas_HarderNarasinhanFilTangent}:
	\begin{Proposition}\emph{(= Proposition \ref{flatnessOnCYs}, c.f.\ \cite[Remark 3.7.(ii)]{biswas_HarderNarasinhanFilTangent})}
	    Let $Y$ be a log terminal projective variety with numerically trivial canonical class $\cc_1(Y)=0$. Then any holomorphic vector bundle on $Y$ admitting a holomorphic connection is necessarily numerically flat.\label{IntroBiswas}
	\end{Proposition}
    As a second application of Proposition \ref{IntroBiswas}, in \cref{sec:SplittingTangentSequence} we give a criterion for a fibration over a $K$-trivial variety to be locally constant which immediately yields the following strengthening of \cite[Theorem 1.1]{hosonoIwaiMatsumura_PsefTangentBundle}:
	\begin{Theorem}\emph{(c.f.\ \cite{hosonoIwaiMatsumura_PsefTangentBundle})} Let $X$ be a smooth projective variety and assume that the tangent bundle $\mathcal{T}_X$ of $X$ admits a positively curved singular Hermitian metric. Then the Albanese $\alpha\colon X\rightarrow \Alb(X)$ is a locally constant fibration. \label{IntroStructureTheoryPsefTangent}
	\end{Theorem}
    Note that in \cite[Theorem 1.1]{hosonoIwaiMatsumura_PsefTangentBundle} it was already proved that $\alpha$ is smooth and that the fibres of $\alpha$ are mutually isomorphic.
	
	The structure of this paper is as follows: In \cref{sec:Preliminaries} we collect some preliminary results on semistability of vector bundles and holomorophic connections in a principal bundle. In particular, we prove Proposition \ref{IntroBiswas}. 
	
	Now, let $f\colon X\rightarrow Y$ be a locally constant fibration as above. In \cref{sec:LocallyConstantFibProjective} we study when $X$ is a projective variety and not merely a complex analytic space in terms of the transition functions of $f$. Then, in \cref{sec:SplittingTangentSequence}, we prove Theorem \ref{IntroStructureTheoryPsefTangent}. Finally, the role of 
    \cref{sec:LocConstFibrAndPositiveCurvature} is to complete the proof of Theorem \ref{IntroCharacterisationMfdsNefACBundle}.

    The original motivation for this work was to construct a counterexample to the Non-Vanishing problem for varieties with nef anti-canonical divisor \cite{lazic_NonVanishingACBundle}; this issue was investigated further in \cite{Mul_NonVanishingACBundle}.
	
	Throughout this text, we work over the field $\C$ of complex numbers and we freely use the standard notions and terminology regarding singularities of algebraic varieties detailed in \cite[Chapter 2]{kollar_BirationalGeometry}.

	\subsection*{Acknowledgements}
	The author would like to express his sincere gratitude towards his advisor, Daniel Greb, for sharing his vast knowledge and his valuable advice. He would also like to warmly thank Shin-Ichi Matsumura for his great talk at the Workshop on Complex Analysis and Geometry in Essen and the helpful discussions afterwards which sparked the author's interest in further investigating this subject. Last but not least many thanks to the anonymous referee whose detailed suggestions helped significantly to clean up the exposition.
	
	While writing this article, the author was partially supported by the DFG research training group 2553 \emph{Symmetries and classifying spaces: analytic, arithmetic and derived}.

	\section{Preliminary results}
	\label{sec:Preliminaries}
	
	\subsection{Selected facts on Semistability and Flatness}
	
	In this subsection, we collecting some basic results concerning semistability of coherent sheaves and prove Proposition \ref{IntroBiswas}. 
 
 Throughout this subsection we fix a normal, complex projective variety $X$ of dimension $n$ and and an ample divisor $H$ on $X$. Following \cite{huybrechts_ComplexGeometry, kobayashi_DiffGeo} we define:
	
	\begin{Definition}
		A torsion-free coherent $\mathcal{O}_X$-module $\mathcal{E}$ is called \emph{$H$-semistable} if 
            $$ \mu(\mathcal{F}) \leq \mu(\mathcal{E}) $$
        for any coherent subsheaf $\mathcal{F}\subseteq \mathcal{E}$ of rank $\rk(\mathcal{F})>0$. Here, the \emph{slope} $\mu(\mathcal{F})$ is defined by $\mu(\mathcal{F})= \frac{\deg(\mathcal{F})}{\rk(\mathcal{F})}$ with $\deg(F) = \cc_1(\mathcal{F})\cdot H^{n-1}$, c.f.\ \emph{\cite[Definition 1.2.11]{huybrechtsLehn_ModuliOfSheaves}}.
        \label{def:semistability}
	\end{Definition}
	\begin{Remark}
		In Definition \ref{def:semistability} it is equivalent to require $\mu(\mathcal{F}) \leq \mu(\mathcal{E})$ only for \emph{saturated} subsheaves $\mathcal{F}\subseteq \mathcal{E}$, i.e.\ such that the quotient $\mathcal{E}/\mathcal{F}$ is torsion-free. Hence, $\mathcal{E}$ is $H$-semistable if and only if $\mu(\mathcal{Q}) \geq \mu(\mathcal{E})$ for any torsion-free quotient sheaf $\mathcal{E}\twoheadrightarrow \mathcal{Q}$ with $\rk\mathcal{Q} > 0$, see \cite[Proposition 5.7.6]{kobayashi_DiffGeoOfVectorBundles}.
	\end{Remark}
	The following two results are classical:
	\begin{Theorem}\emph{(Metha-Ramanathan, Flenner \cite[Theorem 7.1.1]{huybrechtsLehn_ModuliOfSheaves})}
		
		\noindent
		Let $X$ be a normal projective variety and let $\mathcal{E}$ be a torsion-free coherent $\mathcal{O}_X$-module. Then $\mathcal{E}$ is $H$-semistable if and only if $\mathcal{E}|_C$ is semistable for some complete intersection curve $C = H_1\cap\ldots\cap H_{n-1}$ with $H_1, \dots H_{n-1}\in |aH|$ general and  $a\gg0$.\label{mehtaRamanatan}
	\end{Theorem}
    Note that a curve $C$ as in Theorem \ref{mehtaRamanatan} is smooth for dimension reasons and by Bertini's theorem.
    \begin{Theorem}\emph{(Harder-Narasimhan filtration \cite[Theorem 5.7.15]{kobayashi_DiffGeo})}
        
		\noindent
		Let $X$ be a normal projective variety and let $\mathcal{E}$ be a torsion-free coherent $\mathcal{O}_X$-module. Then there exists a unique filtration 
        $$0=\mathcal{E}_0\subsetneq \mathcal{E}_1\subsetneq \ldots \subsetneq \mathcal{E}_{\ell} = \mathcal{E}$$
        by coherent $\mathcal{O}_X$-modules such that the successive quotients $\mathcal{E}_i/\mathcal{E}_{i-1}$ are semistable, torsion-free $\mathcal{O}_X$-modules and the slopes strictly decrease: 
        $$\mu(\mathcal{E}_{i+1}/\mathcal{E}_i) < \mu(\mathcal{E}_i/\mathcal{E}_{i-1}), \qquad \forall i = 1, \ldots, \ell-1.$$
        \label{HNfiltration}
    \end{Theorem}
	We will require the following generalisation of \cite[Remark 3.7.(ii)]{biswas_HarderNarasinhanFilTangent} to the case of normal varieties; the proof is essentially identical but we recall the argument for the convenience of the reader:
	\begin{Lemma}\emph{(Biswas \cite[Remark 3.7.(ii)]{biswas_HarderNarasinhanFilTangent})}
		
		\noindent
		Let $X$ be a normal projective variety. Assume that the tangent sheaf $\mathcal{T}_X := (\Omega_X^1)^*$ of $X$ is $H$-semistable of slope $\mu(\mathcal{T}_X)\geq 0$. Let $\mathcal{E}$ be a vector bundle on $X$ such that $\mathcal{E}|_{X_{\mathrm{reg}}}$ admits a holomorphic connection $D\colon \mathcal{E}|_{X_{\mathrm{reg}}}\rightarrow \mathcal{E}|_{X_{\mathrm{reg}}}\otimes \Omega^1_{X_{\mathrm{reg}}}$. Then $\mathcal{E}$ is $H$-semistable.\label{biswasConnectionsImppliesSemiStable}
	\end{Lemma}
	\begin{proof}
		Let $C = H_1\cap\ldots\cap H_{n-1}$ be a general complete intersection curve for $H$ of sufficiently large degree as in Theorem \ref{mehtaRamanatan}. We will show that $\mathcal{E}|_C$ is semistable. Let
		\begin{align*}
		0=\mathcal{E}_0\subsetneq \mathcal{E}_1\subsetneq \ldots \subsetneq \mathcal{E}_{\ell} = \mathcal{E}|_C
		\end{align*}
		be the Harder-Narasimhan filtration of $\mathcal{E}|_C$ as in Theorem \ref{HNfiltration}. We consider the composition
		\begin{align}
		\mathcal{E}_1 \hookrightarrow \mathcal{E}_{\ell} \overset{D}{\longrightarrow} \mathcal{E}_{\ell} \otimes \Omega_X^1|_C \rightarrow \mathcal{E}_{\ell}/\mathcal{E}_{\ell -1} \otimes \Omega_X^1|_C.
        \label{eq:Biswas1}
		\end{align}
        Note that $C \subseteq X_{\mathrm{reg}}$ by Bertini's theorem so that (\ref{eq:Biswas1}) is well-defined.
		By the Leibniz-rule (\ref{eq:Biswas1}) is $\mathcal{O}_C$-linear. Taking duals yields a map of $\mathcal{O}_C$-modules
		\begin{align*}
		\delta \colon \mathcal{E}_1\otimes \mathcal{T}_X|_C \rightarrow \mathcal{E}_{\ell}/\mathcal{E}_{\ell - 1}.
		\end{align*}
		Now, $\mathcal{T}_X|_C$ is semistable by Theorem \ref{mehtaRamanatan} and $\mathcal{E}_1$ is semistable by assumption. Consequently, by \cite[Theorem 3.1.4]{huybrechtsLehn_ModuliOfSheaves}, also $\mathcal{E}_1\otimes \mathcal{T}_X|_C$ is semistable of slope
		\begin{align*}
		\mu\left( \mathcal{E}_1\otimes \mathcal{T}_X|_C \right) = \mu\left( \mathcal{E}_1 \right) + \mu\left( \mathcal{T}_X|_C \right) = \mu\left( \mathcal{E}_1 \right) + \mu\left( \mathcal{T}_X \right) \geq \mu\left( \mathcal{E}_1 \right).
		\end{align*}
		Since $\textmd{Im}(\delta)$ is a quotient of $\mathcal{E}_1\otimes \mathcal{T}_X|_C$ we deduce that
		\begin{align*}
			\mu\left( \textmd{Im}(\delta) \right) \geq \mu\left(\mathcal{E}_1\otimes \mathcal{T}_X|_C \right) \geq \mu\left( \mathcal{E}_1 \right).
		\end{align*}
		On the other hand, as $\textmd{Im}(\delta) \subseteq \mathcal{E}_{\ell}/\mathcal{E}_{\ell - 1}$ we have that either
		\begin{align*}
		\mu\left(\textmd{Im}(\delta)\right) \leq \mu\left( \mathcal{E}_{\ell}/\mathcal{E}_{\ell - 1} \right) < \ldots < \mu\left(\mathcal{E}_1/\mathcal{E}_0\right) = \mu\left(\mathcal{E}_1\right)
		\end{align*}
		or $\textmd{Im}(\delta)$ is trivial. In effect, $0= \delta$ must be trivial and so $D$ maps $\mathcal{E}_1$ into $\mathcal{E}_{\ell - 1} \otimes \Omega_X^1|_C$. Repeating the above argument, we find an $\mathcal{O}_C$-linear map 
        $$\delta'\colon \mathcal{E}_1\otimes \mathcal{T}_X|_C \rightarrow \mathcal{E}_{\ell -1}/\mathcal{E}_{\ell - 2}$$
        and, as before, we see that $\delta' = 0$. Arguing inductively, we find that $D$ restricts to a connection $D_1\colon \mathcal{E}_1\rightarrow \mathcal{E}_1\otimes \Omega^1_X|_C\rightarrow \mathcal{E}_1\otimes \Omega^1_C$ on the semistable vector bundle $\mathcal{E}_1$. As all Chern classes of a bundle with a holomorphic connection vanish \cite[Theorem 4]{atiyah_ConnectionsInFibreBundles} we see that $\mu(\mathcal{E}_1) = \mu(\mathcal{E}|_C) = 0$. Finally, applying the above argument to the induced map 
        $$D'\colon \left(\mathcal{E}|_C/\mathcal{E}_1 \right) \rightarrow \left(\mathcal{E}|_C/\mathcal{E}_1\right) \otimes \Omega_X^1|_C$$
        we find that $D'$ restricts to a holomorphic connection on $\mathcal{E}_2$ and, hence, that $\mu(\mathcal{E}_2) = 0$. By Theorem \ref{HNfiltration} it follows that $\mathcal{E}_1 = \mathcal{E}|_C$ is semistable. An application of Theorem \ref{mehtaRamanatan} now yields the semistability of $\mathcal{E}$.
	\end{proof}
	The reason we are interested in semistability is the following result which allows to infer flatness of bundles. It relies crucially on the recent extension of the \emph{Non-Abelian Hodge Correspondence} to klt varieties due to Greb-Kebekus-Peternell-Taji.
	\begin{Theorem}\emph{(c.f.\ \cite{gkpt_HodgeTheoremKltSpaces})}
		Let $(X, \Delta)$ be a projective klt pair of dimension $n$ and fix an ample divisor $H$ on $X$. For any vector bundle $\mathcal{E}$ on $X$ the following are equivalent:
		\begin{itemize}
			\item[(1)] The bundles $\mathcal{E}$ and $\mathcal{E}^*$ are both nef,
			\item[(2)] $\mathcal{E}$ is $H$-semistable and holomorphically flat, i.e.\ there exists a $\C$-local system $\E$ on $X$ such that $\mathcal{E} \cong \E\otimes_{\C} \mathcal{O}_X$ as holomorphic vector bundles,
			\item[(3)] $\mathcal{E}$ is $H$-semistable and satisfies $\cc_1(\mathcal{E})\cdot H^{n-1} = \ch_2(\mathcal{E})\cdot H^{n-2} = 0$.
		\end{itemize}
		In case $\mathcal{E}$ satisfies one of these properties we will say that $\mathcal{E}$ is \emph{numerically flat}.\label{druelNumFlat}
	\end{Theorem}
	\begin{proof}
		The equivalence of $(1)$ and $(3)$ is readily deduced from the main statement of \cite[Theorem 3.4]{gkpt_HodgeTheoremKltSpaces} and \cite[Theorem 1.18]{DPS_ManifoldsWithNefTangentBundle}, see \cite[Theorem 2.9]{druel_NumFlatBundlesKltSpaces} for the precise argument. Moreover, clearly $(2) \Rightarrow (3)$.

        Finally, the converse implication $(3) \Rightarrow (2)$ is a direct consequence of \cite[Theorem 3.4]{gkpt_HodgeTheoremKltSpaces} and \cite[Section 3]{simpson_NAHT}. Indeed, fix a resolution $\pi\colon \widehat{X}\rightarrow X$. According to \cite[Theorem 3.4, Remark 3.7]{gkpt_HodgeTheoremKltSpaces}, there exists a local system $\E$ on $X$ such that $(\pi^*\mathcal{E}, 0)$ is a semistable Higgs bundle with vanishing Chern classes and whose underlying $C^{\infty}$-vector bundle $E$ is the same as the underlying $C^{\infty}$-bundle of $\pi^*(\E\otimes_{\C} \mathcal{O}_X)$. Now, a priori these two holomorphic structures on $E$ might differ but by the remark after \cite[Lemma 3.5]{simpson_NAHT} they do not, precisely because the Higgs field vanishes. In conclusion, $\pi^*\mathcal{E} \cong \pi^*(\E\otimes_{\C} \mathcal{O}_X)$ and it follows by the projection formula that $\mathcal{E} \cong \E\otimes_{\C} \mathcal{O}_X$ as proclaimed.
	\end{proof}
    Note that a holomorphic vector bundle on a smooth variety is holomorphically flat if and only if it admits a flat holomorphic connection, c.f.\ Remark \ref{rem:FlatPrincipalBundles}.
    
	The following corollary ( = Proposition \ref{IntroBiswas}) will be crucial in the later sections:
	
	\begin{Proposition}
		Let $Y$ be a log terminal projective variety with numerically trivial canonical class $\cc_1(Y)=0$. Let $\mathcal{E}$ be a vector bundle on $Y$ and assume that $\mathcal{E}|_{Y_{\mathrm{reg}}}$ admits a holomorphic connection.
		\begin{itemize}
			\item[(1)] The bundle $\mathcal{E}$ is semistable with respect to any ample divisor on $Y$.
			\item[(2)] In fact, $\mathcal{E}$ is numerically flat.
		\end{itemize}
	In particular, a vector bundle on $Y$ is numerically flat if and only if its restriction to $Y_{\mathrm{reg}}$ is holomorphically flat.\label{flatnessOnCYs}
	\end{Proposition}
	\begin{proof}
		According to \cite{guenancia_SemistabilityOfTangent} it holds that $\mathcal{T}_Y$ is semistable of slope $\mu(\mathcal{T}_Y) = 0$. Thus, item $(1)$ immediately follows from Lemma \ref{biswasConnectionsImppliesSemiStable}. Regarding $(2)$, by \cite[Theorem 1.14]{GKP_MaximalllyQuasiEtaleCovers} there exists a finite quasi-\'etale cover $\pi\colon Y'\rightarrow Y$, of degree $d$ say, such that $\pi^*\mathcal{E}$ is holomorphically flat. Using $(1)$, Theorem \ref{druelNumFlat} shows that $\pi^*\mathcal{E}$ is numerically flat. In particular, $\cc_1(\mathcal{E})\cdot H^{n-1} = \frac{1}{d}\cdot \cc_1(\pi^*\mathcal{E})\cdot (\pi^*H)^{n-1} = 0$ and $\ch_2(\mathcal{E})\cdot H^{n-2} = \frac{1}{d}\cdot \ch_2(\pi^*\mathcal{E})\cdot (\pi^*H)^{n-2} = 0$. Thus, $\mathcal{E}$ is numerically flat by Theorem \ref{druelNumFlat}. The addendum is an immediate consequence of Lemma \ref{biswasConnectionsImppliesSemiStable} and Theorem \ref{druelNumFlat}.
	\end{proof}

	\subsection{Some words on Principal Bundles}

    \begin{Definition}
        Let $F$ be an (irreducible) complex analytic variety. A \emph{holomorphic fibre bundle with fibre $F$} (also known as \emph{locally trivial fibration} in the literature) is a holomorphic map $f\colon X\rightarrow Y$ between complex analytic varieties such that there exists an open cover $Y = \bigcup U_i$ in the Euclidean topology and isomorphims
        \begin{align*}
            \theta_i\colon U_i \times F \rightarrow X|_{U_i} := f^{-1}(U_i)
        \end{align*}
        over $U_i$. In this case, all fibres of $f$ are isomorphic to $F$. In particular, there exist uniquely determined maps $g_{i,j}\colon U_i\cap U_j\rightarrow \Aut(F)$ such that $\theta_i(\theta_j^{-1}(x, y)) = (g_{i, j}(y)(x), y)$ for all $x\in F,\ y\in U_i\cap U_j$. We call $(U_i, \theta_i)$ a \emph{trivialising open cover} for $X\rightarrow Y$ with \emph{transition functions} $(g_{i,j})$. Holomorphic fibre bundles are also known as \emph{locally trivial fibrations} in the literature.
    \end{Definition}
    In case $F$ is compact the group of holomorphic automorphisms $G:= \Aut(F)$ is a complex Lie group \cite[Section 2.3]{akhiezer_AutsOfComplexMfds}. Moreover, the transition functions $g_{i,j}\colon U_i\cap U_j\rightarrow \Aut(F)$ are holomorphic.
    
    Now, let $G\subseteq \Aut(F)$ be a subgroup. We say the \emph{structure group} of $f\colon X\rightarrow Y$ \emph{can be reduced to} $G$ if we can choose the $(U_i, \theta_i)$ in such a way that the $g_{i,j}$ take values in $G$. In this case the $(g_{i,j})$ form a \v{C}ech cocyle in $\check{\mathrm{H}}^1(Y,G)$, which is independent of the choice of $(U_i, \theta_i)$. Conversely, any \v{C}ech cocyle $(g_{i,j}) \in \check{\mathrm{H}}^1(Y,G)$ can be used to glue the trivial fibre bundles $U_i\times F \rightarrow U_i$ into a fibre bundle $f\colon X\rightarrow Y$ with fibre $F$, structure group $G$ and transition functions $(g_{i,j})$. It is well-known that a fibre bundle is determined uniquely by its corresponding class in $\check{\mathrm{H}}^1(Y,G)$.
     \begin{Definition}
		A proper map $f\colon X \rightarrow Y$ between complex analytic varieties is a \emph{locally constant fibre bundle} with fibre $F$ if there exists a group homomorphism $\rho\colon\pi_1(Y) \rightarrow \Aut(F)$ and an isomorphism of complex analytic varieties
		\begin{align*}
		X \cong (\widetilde{Y}\times F)/\pi_1(Y)
		\end{align*}
		over $Y$, where $\widetilde{Y}\rightarrow Y$ denotes the universal cover of $Y$. Here, $\pi_1(Y)$ acts on $\widetilde{Y}$ in the natural way and on $F$ via $\rho$. Moreover, $f$ is called \emph{locally constant} with respect to some $\mathds{Q}$-Weil divisor $\Delta$ on $X$ if
	   \begin{itemize}
		  \item[(i)] $\Delta$ is \emph{horizontal}, i.e.\ any irreducible component of $\Delta$ is dominant over $Y$. Equivalently, $\Delta_F := \Delta\cap F$ is a divisor on $F$.
		  \item[(ii)] There exists a group homomorphism $\rho\colon\pi_1(Y) \rightarrow \Aut(F, \Delta_F) := \{ \varphi \in \Aut(F)  | \ \varphi(\Delta_F) = \Delta_F \}$ such that $(X, \Delta) \cong (\widetilde{Y}\times F,\ pr_F^*\Delta_F)/\pi_1(Y)$.
	   \end{itemize}
		
		In other words, $f$ is locally constant if and only if it is a holomorphic fibre bundle whose transition functions can be chosen to be locally constant. Note that locally constant fibre bundles are also known as \emph{flat} fibre bundles in the literature.\label{ex:LocallyConstantFibration}
	\end{Definition}
    \begin{Example}
       Let $F = G$ be a complex Lie group. Then $G$ acts on itself by left multiplication and so we can view $G$ as a subgroup of $\Aut(G)$. A holomorphic fibre bundle $p\colon \mathcal{G}\rightarrow Y$ with fibre $G$ and structure group $G\subseteq \Aut(G)$ is also called a \emph{holomorphic principal G-bundle}. In this case, it holds that $Y \cong \mathcal{G}/G$. Principal bundles will play a crucial role in the later sections; for us, they will all arise in the following way:
    \end{Example}
	\begin{Construction}
		Let $f\colon X\rightarrow Y$ be a proper holomorphic fibre bundle with fibre $F$. We denote by $\mathcal{G}(f)\rightarrow Y$ the unique holomorphic principal $G := \Aut(F)$-bundle defined using the same transition functions as $f$. More explicitly, let $Y = \bigcup U_i$ be any trivialising open cover for $f$ with transition functions $g_{i,j}\colon U_i\cap U_j\rightarrow G$. Then, $(g_{i,j})$ form a \v{C}ech cocyle in $\check{\mathrm{H}}^1(Y,G)$ which gives rise to $\mathcal{G}(f)$. Note that $\mathcal{G}(f)$ does not depend (up to isomorphism) on the choice of $(g_{i,j})$.
  
        Conversely, given a complex Lie group $G$, a holomorphic principal $G$-bundle $p\colon \mathcal{G}\rightarrow Y$ and an effective action of $G$ on another complex analytic space $F$ we can consider $G\subseteq \Aut(F)$ as a subgroup. Then any set of transition functions $g_{i,j}$ for $\mathcal{G}$ can be used to define a holomorphic fibre bundle $f\colon \mathcal{G}\times^G F\rightarrow Y$ with fibre $F$. Again, this bundle does not depend on the choice of $g_{i,j}$. We call $\mathcal{G}\times^G F\rightarrow Y$ the bundle \emph{associated to $\mathcal{G}$ with fibre $F$}.
        
        For example, in the notation above it holds that $X = \mathcal{G}(f)\times^G F$.\label{con:AssociatedPrincipalBundle}
	\end{Construction}
	
	\begin{Remark}(Holomorphic Connections in Principal Bundles)
		
		\noindent
		Following \cite{atiyah_ConnectionsInFibreBundles}, a holomorphic connection in a principal $G$-bundle $p\colon \mathcal{G}\rightarrow Y$ over a smooth variety $Y$ is nothing but a holomorphic (!) splitting of the short exact sequence
		\begin{align*}
		0 \rightarrow \mathfrak{ad}(\mathcal{G}) \rightarrow \left(p_*\mathcal{T}_{\mathcal{G}}\right)^G \rightarrow \mathcal{T}_Y \rightarrow 0.
		\end{align*}
		Here, $G$ is an arbitrary complex Lie group, $\mathfrak{ad}(\mathcal{G}):= \mathcal{G}\times^G \mathfrak{g}$ is the vector bundle associated to the \emph{adjoint representation} $G\rightarrow \GL(\mathfrak{g})$ and $\left(p_*\mathcal{T}_{\mathcal{G}}\right)^G \subseteq p_*\mathcal{T}_{\mathcal{G}}$ is the subbundle of sections which are invariant under the natural action of $G$. Equivalently, a holomorphic connection is a holomorphic $G$-invariant subbundle $\mathcal{T}\subseteq \mathcal{T}_{\mathcal{G}}$ such that $\mathcal{T}_{\mathcal{G}} = \mathcal{T}\oplus \mathcal{T}_{\mathcal{G}/Y}$. The equivalence follows from the fact that $\mathfrak{ad}(\mathcal{G}) = (p_*\mathcal{T}_{\mathcal{G}/Y})^G$ which is straightforward to verify.
		
		Given a holomorphic connection in $\mathcal{G}$ one defines the parallel transport maps in the usual way, see for example \cite[Section II.3]{kobayashi_DiffGeo}. We say that the connection is flat if and only if the parallel transport maps are invariant under homotopy. Then $\mathcal{G}$ admits a flat holomorphic connection if and only if it is flat in the sense that there exists a group homomorphism $\rho\colon\pi_1(Y)\rightarrow G$ such that $\mathcal{G} = (\widetilde{Y}\times G)/\pi_1(Y)$. Here, $\widetilde{Y}$ denotes the universal cover of $Y$ and $\pi_1(Y)$ acts on $\widetilde{Y}$ in the natural way and on $G$ through $\rho$. Equivalently, $\mathcal{G}$ is holomorphically flat if and only if it can be defined using locally constant transition functions. Compare also to the treatment in \cite[Proposition 14]{atiyah_ConnectionsInFibreBundles} for proofs.\label{rem:FlatPrincipalBundles}
	\end{Remark}

	\section{Projectivity of Locally Constant Fibre Bundles}
	\label{sec:LocallyConstantFibProjective}

	\noindent
	This section is concerned with the projectivity statement in Theorem \ref{IntroCharacterisationMfdsNefACBundle}. This material seems to be known to the experts in the field but our statements are slightly more optimal then what one can usually find in the literature.
	\begin{Definition}
		Let $X$ be a complex analytic space and let $G\subseteq \Aut(X)$ be a subgroup. A line bundle $\mathscr{L}$ on $X$ is \emph{$G$-invariant} if $g^*\mathscr{L}\cong \mathscr{L}$ for any $g\in G$. In this case we say that $\mathscr{L}$ is \emph{$G$-linearisable} if the action of $G$ on $X$ can be lifted to an action of $G$ on the total space of $\mathscr{L}$ via bundle automorphisms.
	\end{Definition}
	\begin{Example}
		Let $V$ be a $\C$-vector space of dimension $n$. Then the line bundle $\mathcal{O}_{\mP(V)}(1)$ is clearly invariant under the full automorphism group $\Aut(\mP(V)) = \mathrm{PGL}(V)$ but this action is \emph{not} linearisable, see for example \cite[Example 4.2.4]{brion_LinearisationsOfGroupActions}. 
  
        On the other hand, $\mathcal{O}_{\mP(V)}(n)$ \emph{does} admit a natural $\mathrm{PGL}(V)$-linearisation: Indeed, $\mathcal{O}_{\mP(V)}(n) = \mathcal{O}_{\mP(V)}(-K_{\mP(V)})$ and so $\varphi\in \Aut(\mP(V))$ acts on $\mathcal{O}_{\mP(V)}(n)$ via $\det(d\varphi)$.\label{ex:LinearisationsLineBundles}
	\end{Example}
	The following result shows that this picture generalises to any projective variety; it is usually only stated in case $\mathrm{Stab}(X, \mathscr{L})$ is connected \cite{brion_LinearisationsOfGroupActions, mumford_GIT}:
	\begin{Lemma}
		Let $X$ be a normal projective variety and let $\mathscr{L}$ be an ample line bundle on $X$. Then the stabiliser
		\begin{align*}
		\mathrm{Stab}(X, \mathscr{L}) := \left\{ g\in \Aut(X) |\ g^*\mathscr{L}\cong \mathscr{L} \right\} \subseteq \Aut(X)
		\end{align*}
		is a linear algebraic group. Moreover, there exists some integer $m>0$ such that $\mathscr{L}^{\otimes m}$ is $\mathrm{Stab}(X, \mathscr{L})$-linearisable. \label{linearisationsAmpleBundles}
	\end{Lemma}
	\begin{proof}
		Fix an integer $k>0$ for which $\mathscr{L}^{\otimes k}$ is very ample. We abbreviate $V:= \mathrm{H}^0(X, \mathscr{L}^{\otimes k})$, $n:=\dim V$ and $G:= \mathrm{Stab}\left(X, \mathscr{L}^{\otimes k}\right)$. Since $X$ is compact we have $\Aut_{\mathcal{O}_{X}}(\mathscr{L}) = \mathcal{O}_X^\times(X) = \mathds{C}^{\times}$ and we infer that 
		\begin{align*}
			\rho\colon & G \rightarrow \mathrm{PGL}(V) = \Aut(\mP(V)),\\
            & g\mapsto \left[g^*\colon \mathrm{H}^0\Big(X, \mathscr{L}^{\otimes k}\Big)\rightarrow \mathrm{H}^0\Big(X, g^*\mathscr{L}^{\otimes k}\Big) \cong \mathrm{H}^0\Big(X, \mathscr{L}^{\otimes k}\Big) \right]
		\end{align*}
		 is a well-defined group homomorphism, independent of choices of isomorphisms $ g^*\mathscr{L}^{\otimes k} \cong  \mathscr{L}^{\otimes k}$ for $g\in G$. By construction the embedding $\iota\colon X \hookrightarrow \mP(V)$ is $\rho$-equivariant and it follows that
		\begin{align*}
		    \rho\colon  \mathrm{Stab}\Big(X, \mathscr{L}^{\otimes k}\Big) = G \hookrightarrow \Aut(\mP(V))
		\end{align*}
		identifies $G$ with the subgroup of $\mathrm{PGL}(V)$ leaving $X\hookrightarrow \mP(V)$ invariant. We infer that $G$ is linear algebraic: Indeed, if $X$ is cut out in $\mP(V)$ by the homogeneous polynomials $f_1,\ldots, f_s$, then
		\begin{align*}
		    G = \bigcap_{x\in X} \Big\{ g\in \mathrm{PGL}(V) |\ f_1(gx)= \ldots = f_s(gx) = 0 \Big\} \subseteq \mathrm{PGL}(V)
		\end{align*}
		and so $G$ is a closed subgroup of the linear algebraic group $\mathrm{PGL}(V)$. As $\mathrm{Stab}(X, \mathscr{L})\subseteq \mathrm{Stab}(X, \mathscr{L}^{\otimes k})$ is a subgroup of finite index - for example because the quotient injects into the finite set of $k$-torsion points of $\mathrm{Pic}(X)$ - if follows that also $\mathrm{Stab} \left(X, \mathscr{L}\right)$ is linear algebraic.
		
		Regarding the second statement, according to Example \ref{ex:LinearisationsLineBundles}, $\mathcal{O}_{\mP(V)}(n)$ admits a natural $\Aut(\mP(V))$-linearisation. It follows that $\mathcal{O}_{\mP(V)}(n)$ is also linearised for the action of the subgroup $\mathrm{Stab}(X, \mathscr{L})\subseteq \Aut(\mP(V))$. Clearly, $\mathscr{L}^{\otimes nk} = \mathcal{O}_{\mP(V)}(n)|_X$ inherits this $\mathrm{Stab}(X, \mathscr{L})$-linearisation and so (letting $m:=nk$) we are done.
	\end{proof}
	The following criterion is essentially elementary:
	\begin{Lemma}
		Let $f\colon X\rightarrow Y$ be a locally constant holomorphic fibre bundle with fibre $F$, given by some group homomorphism $\rho\colon\pi_1(Y) \rightarrow \Aut(F)$. Assume, that $F$ is a normal projective variety and that $Y$ is compact or algebraic. Consider the following assertions:
		\begin{itemize}
			\item[(1)] The morphism $f$ is projective,
			\item[(2)] there exists a $\pi_1(Y)$-invariant ample line bundle on $F$,
			\item[(3)] there exists a $\pi_1(Y)$-invariant and linearisable ample line bundle on $F$,
			\item[(4)] the image of the map $\pi_1(Y) \overset{\rho}{\rightarrow} \Aut(F)/\Aut^0(F)$ is finite,
			\item[(5)] the image of the induced map $\pi_1(Y) \rightarrow \GL(N^1(F))$ is finite.
		\end{itemize}\label{criterionForProjectivity}
        Then $(1) \Leftrightarrow (2) \Leftrightarrow (3) \Rightarrow (4) \Rightarrow (5)$. Moreover, in case $\mathrm{Pic}^0(F) = 0$ then all five assertions are equivalent. Here, $\mathrm{Pic}^0(F)$ denotes the identity component of $\mathrm{Pic}(F)$.
	\end{Lemma}
	Recall that $N^1(F)$, the \emph{N\'eron-Severi group} of $F$, denotes the group of line bundles on $F$ modulo numerical equivalence, or, equivalently, the image of the group homomorphism  $c_1\colon \mathrm{Pic}(F) \rightarrow \mathrm{H}^2(X, \R)$. Moreover, $\Aut^0(F)$ denotes the connected component of the identity in $\Aut(F)$. Note that $\Aut^0(F)\subseteq \Aut(F)$ is always a normal complex Lie subgroup.
	\begin{proof}
	    First of all, if $\mathscr{L}$ is any $f$-ample line bundle on $X$, then we consider the line bundle $p^*\mathscr{L}$ on $X\times_Y \widetilde{Y} \cong \widetilde{Y}\times F\overset{p}{\rightarrow} X$. Note that for any $y\in \widetilde{Y}$ and any $\gamma \in \pi_1(Y, p(y))$ the line bundles $p^*\mathscr{L}|_{\{y\} \times F}$ and $p^*\mathscr{L}|_{\{\gamma \cdot y\} \times F}$ on $F$ can both be canonically identified with $\mathscr{L}|_{F}$. In particular, using that $p\circ \gamma = p$, we find
        $$
        \gamma^*\mathscr{L}|_F
        \cong \gamma^*\Big(p^*\mathscr{L}|_{\{\gamma \cdot y\} \times F}\Big)
        = \big(\gamma^*p^*\mathscr{L}\big)|_{\{y\} \times F}
        = (p^*\mathscr{L})|_{\{y\} \times F}
        \cong \gamma^*\mathscr{L}|_F.
        $$
        In other words, $\mathscr{L}|_{F}$ is $\pi_1(Y)$-invariant. This proves that $(1)\Rightarrow (2)$. 
	    
	    The implication $(2) \rightarrow (3)$ is clear in view of Lemma \ref{linearisationsAmpleBundles}. Moreover, if $\mathscr{L}_F$ is a $\pi_1(Y)$-linearisable ample line bundle on $F$ then $pr_F^*\mathscr{L}_F$ is a $\pi_1(Y)$-invariant linearisable line bundle on $\widetilde{Y}\times F$ and, thus, its quotient by the action of $\pi_1(Y)$ is an $f$-ample holomorphic line bundle on $X$. Hence, $f$ is projective. This finishes the proof that $(1), (2)$ and $(3)$ are pairwise equivalent.
        
        Regarding $(3)\Rightarrow (4)$, let $\mathscr{L}_F$ be a $\pi_1(Y)$-invariant ample line bundle on $F$. Then the stabiliser $\mathrm{Stab}(F, \mathscr{L}_F)$ is an algebraic group by Lemma \ref{linearisationsAmpleBundles} and, thus, the group
		\begin{align*}
		\mathrm{Stab}\left(F, \mathscr{L}_F\right)/ \mathrm{Stab}^0\left(F, \mathscr{L}_F\right)
		\end{align*}
		is finite. Since $\rho\colon\pi_1(Y) \rightarrow \Aut(F)$ has image in $\mathrm{Stab}\left(F, \mathscr{L}_F\right)$ by assumption we see that $\pi_1(Y) \rightarrow \Aut(F)/\Aut^0(F)$ factors through $\mathrm{Stab}\left(F, \mathscr{L}_F\right)/\mathrm{Stab}^0\left(F, \mathscr{L}_F\right)$. We conclude that $\pi_1(Y) \overset{\rho}{\rightarrow} \Aut(F)/\Aut^0(F)$ has finite image, thereby confirming $(4)$. The implication $(4)\Rightarrow (5)$ is tautologous as $\pi_1(Y) \rightarrow \GL(N^1(F))$ factors through $\Aut(F)/\Aut^0(F)$. 
  
        Finally, assuming that $\mathrm{Pic}^0(F) = 0$, let us show that $(5)\Rightarrow (2)$: Fix an ample line bundle $\mathscr{L}'_F$ on $F$ and let us denote $H:= \textmd{Im}\big(\pi_1(Y) \rightarrow \GL(N^1(F))\big)$. Then the numerical class
		\begin{align*}
		\sum_{\gamma\in H} \Big[\gamma^*\mathscr{L}'_F\Big] \in N^1(F)\subseteq \mathrm{H}^2(X, \R)
		\end{align*}
		belongs to an ample line bundle $\mathscr{L}_F$ whose numerical class is $\pi_1(Y)$-invariant. Replacing $\mathscr{L}_F$ by some multiple (to get rid of the torsion) if necessary we can assume that 
		\begin{align*}
		    [\mathscr{L}_F]\in \mathrm{Pic}(F) = \mathrm{Pic}(F)/\mathrm{Pic}^0(F)\hookrightarrow \mathrm{H}^2(X, \Z)
		\end{align*}
		is itself invariant, finishing the proof.
	\end{proof}
	\begin{Corollary}
	    Let $Y$ be a projective log terminal variety with numerically trivial canonical class $\cc_1(Y) = 0$, let $(F, \Delta_F)$ be a projective klt pair with $F$ being rationally connected and let $\rho\colon \pi_1(Y)\rightarrow \Aut(F, \Delta_F)$ be a group homomorphism. Write $(X, \Delta) := (\widetilde{Y}\times F,\ pr_F^*\Delta_F)$ and let $f\colon X\rightarrow Y$ denote the corresponding locally constant fibration.
		
		Then $X$ is projective if and only if the image of $\pi_1(Y)\rightarrow \Aut(F)/\Aut^0(F)$ is finite. In particular, this is true provided one of the following assertions is satisfied:
		\begin{itemize}
			\item[(1)] The anti-log canonical divisor $-(K_F+\Delta_F)$ is big.
			\item[(2)] The augmented irregularity $\widetilde{q}(Y) = 0$ vanishes. Equivalently, the singular Beauville-Bogomolov decomposition of $Y$ does not contain a torus factor.
		\end{itemize}\label{criterionForProjectivitySpecialCase}
	\end{Corollary}
    Here, $\widetilde{q}(Y) := \sup\{ \ h^1(Y', \mathcal{O}_{Y'})\ | \ \nu\colon Y'\rightarrow Y \textmd{ finite quasi-\'etale} \}$, see \cite[Definition 2.20]{grebGuenanicaKebekus_KltCalabiYaus}.
	\begin{proof}
		Note that under our assumptions the singularities of $F$ are rational \cite[Theorem 5.22]{kollar_BirationalGeometry}. Thus, $\dim \textmd{Pic}^0(F) = q(F) = 0$, $F$ being rationally connected, and Lemma \ref{criterionForProjectivity} applies. This proves the first statement.
		
		Now, assuming that $-(K_F+\Delta_F)$ is big the group $\Aut(F, \Delta_F)$  is always linear algebraic, see for example the argument given in \cite[Proposition 2.26]{brion_AutomorphismGroups}. In particular, the group $\Aut(F, \Delta|_F)/\Aut^0(F, \Delta|_F)$ is finite and $(1)$ follows.
		
		Regarding the second assertion, if $\widetilde{q}(Y) = 0$ then \emph{any} linear representation of $\pi_1(Y)$ has finite image \cite[Theorem I]{grebGuenanicaKebekus_KltCalabiYaus}. Thus, the image of the group homomorphism $\pi_1(Y)\rightarrow \GL(N^1(F)) \hookrightarrow \GL(N^1(F)_{\Q})$ is finite and we again conclude by using Lemma \ref{criterionForProjectivity}.
	\end{proof}

	\section{Fibrations and Splittings of the relative Tangent Sequence}
	\label{sec:SplittingTangentSequence}
	
	This section is devoted to the proof of the following result; its statement is somewhat surprising to the author since we need not impose any integrability conditions on the subsheaf $f^*\mathcal{T}_Y \subseteq \mathcal{T}_X$:
	\begin{Lemma}
		Let $f\colon X\rightarrow Y$ be a holomorphic fibre bundle with fibre $F$. Assume that $f$ is projective, that $Y$ is a log terminal projective variety with numerically trivial canonical class $\cc_1(Y)=0$ and that $F$ is a normal projective variety such that $\emph{Pic}^0(F) = 0$. If the relative tangent bundle sequence
		\begin{align*}
		0\rightarrow \mathcal{T}_{X/Y} \rightarrow \mathcal{T}_X\rightarrow f^*\mathcal{T}_Y \rightarrow 0
		\end{align*}
		splits holomorphically then $f$ is a locally constant fibration.\label{splittingTangentSequence}
	\end{Lemma}
	\begin{Remark}
	    Note that the converse to Lemma \ref{splittingTangentSequence} also holds true: Given a locally constant fibration $f\colon X= (\widetilde{Y}\times F)/\pi_1(Y)\rightarrow Y$ the relative tangent bundle sequence $0\rightarrow \mathcal{T}_{X/Y} \rightarrow \mathcal{T}_X\rightarrow f^*\mathcal{T}_Y \rightarrow 0$ splits. Indeed, the subsheaf $pr_{\widetilde{Y}}^*\mathcal{T}_{\widetilde{Y}}\subseteq \mathcal{T}_{\widetilde{Y}\times F}$ is clearly $\pi_1(Y)$-invariant and descends to a holomorphic subsheaf of $\mathcal{T}_X$ splitting the relative tangent sequence.
	    \label{rem:SplittingTangentSequence}
	\end{Remark}
	We immedeatly deduce the following strengthening of \cite[Theorem 1.1]{hosonoIwaiMatsumura_PsefTangentBundle} where it was already proved that $\alpha$ is a holomorphic fibre bundle: 
	\begin{Theorem}\emph{(= Theorem \ref{IntroStructureTheoryPsefTangent}, c.f.\ \cite{hosonoIwaiMatsumura_PsefTangentBundle})}
		
		\noindent
		Let $X$ be a smooth projective variety and assume that the tangent bundle $\mathcal{T}_X$ of $X$ admits a positively curved singular Hermitian metric. Then the Albanese morphism $\alpha: X \rightarrow T := \Alb(X)$ is a locally constant fibration.\label{StructureMfdsPsefTangent}
	\end{Theorem}
    \begin{proof}
        By \cite[Theorem 1.1]{hosonoIwaiMatsumura_PsefTangentBundle} there exists a morphism $f\colon X \rightarrow Y$ with connected fibres onto a finite \'etale quotient of an abelian variety $Y$ such that $f$ is a holomorphic fibre bundle. Moreover, the fibres $F$ of $f$ are rationally connected. In particular, by the rigidity lemma \cite[Lemma 1.15]{Debarre_AlgebraicGeometry}, $\alpha$ factors via $\alpha_Y \colon Y \rightarrow T = \Alb(Y)$. According to Lemma \ref{splittingTangentSequence}, $f$ is a locally constant fibration, given by a group homomorphism $\rho\colon \pi_1(Y) \rightarrow \Aut(F)$ say. Moreover, $\alpha_Y$ is a holomorphic fibre bundle, with fibre $F_Y$ say, see for example \cite{cao_NefAnticanonicalBundleII}.

        We claim that Lemma \ref{splittingTangentSequence} applies to $\alpha$. Indeed, $\alpha$ is smooth as the composition of the smooth morphisms $\alpha_Y, f$. Moreover, the fibres $F_X$ of $\alpha$ are mutually isomorphic to the locally constant fibration with base $F_Y$, fibre $F$ and transition functions 
        $$\pi_1(F_Y) \rightarrow \pi_1(Y) \overset{\rho}{\rightarrow} \Aut(F).$$
        Hence, $\alpha$ is a holomorphic fibre bundle \cite{grauertfischer_holomorphicFibreBundles}. Moreover, $\textmd{Pic}^0(F_X) = 0$ by the argument in \cite[Proposition 3.12.(iii)]{DPS_ManifoldsWithNefTangentBundle}. Finally, the short exact sequence
        \begin{align*}
		0\rightarrow \mathcal{T}_{X/T} \rightarrow \mathcal{T}_X\rightarrow f^*\mathcal{T}_T \rightarrow 0
		\end{align*}
        splits holomorphically by \cite[Theorem 1.4]{hosonoIwaiMatsumura_PsefTangentBundle}. Lemma \ref{splittingTangentSequence} concludes the proof.
    \end{proof}
    \begin{Remark}
        According to \cite{hosonoIwaiMatsumura_PsefTangentBundle}, up to replacing $X$ by some finite \'etale cover, it holds that the fibres of $\alpha$ in Theorem \ref{StructureMfdsPsefTangent} are rationally connected. This picture is extremely reminiscent of the results proved by Demailly-Peternell-Schneider and Cao for varieties with nef tangent bundle
        \cite{DPS_ManifoldsWithNefTangentBundle, cao_PhDThesis}. Now, Campana-Peternell \cite{campanaPeternell_Conjecture} famously conjectured that a smooth projective rationally connected variety with nef tangent bundle is homogenoues. In a similar spirit one might wonder how far away a smooth projective rationally connected variety $F$ whose tangent bundle admits a positively curved singular Hermitian metric is from being almost homogeneous? Note that in \cite[Problem 3.12]{hosonoIwaiMatsumura_PsefTangentBundle} the weaker question of whether $-K_F$ must be big in this situation was raised, see \cite[Theorem 1.2]{BaohuaZhang_AlmostHomogeneousVarieties}.
	\end{Remark}
	The proof of Lemma \ref{splittingTangentSequence} will take up the rest of this section. It will be split up into three steps:
	\begin{itemize}
		\item[(1)] Prove the existence of an $f$-ample line bundle $\mathscr{L}$ on $X$ such that $f_*\mathscr{L}$ is associated to $\mathcal{G}(f)$. Here, the principal bundle $\mathcal{G}(f)$ was defined in Construction \ref{con:AssociatedPrincipalBundle}.
		\item[(2)] Translate the existence of a splitting of the relative tangent sequence for $f$ into the existence of a holomorphic connection in $\mathcal{G}(f)|_{Y_{\mathrm{reg}}}$.
		\item[(3)] Conclusion: From $(2)$ we deduce that $f_*\mathscr{L}$ admits a holomorphic connection and, hence, by Proposition \ref{flatnessOnCYs} that $f_*\mathscr{L}$ is numerically flat. This forces $f$ to be locally constant by a standard criterion \cite{cao_PhDThesis}.
	\end{itemize}

	\subsection{Step 1: Existence of a Relatively Ample and Flat Line Bundle}
	
	Throughout this subsection we fix a proper holomorphic fibre bundle $f\colon X\rightarrow Y$ with fibre $F$. Let us abbreviate $G := \Aut(F)$ and $\mathcal{G}:=\mathcal{G}(f)$.
	\begin{Example}
		Suppose that $F$ admits an ample $G = \Aut(F)$-linearisable line bundle $\mathscr{L}_F$. Then $\mathscr{L}:= \mathcal{G}\times^G \mathscr{L}_F$ makes sense as a holomorphic line bundle on $X$. Moreover, $\mathscr{L}$ is clearly $f$-ample and $f_*\mathscr{L} = \mathcal{G}\times^G \mathrm{H}^0(F, \mathcal{O}_F(-K_F))$ is associated to $f$. This construction applies for example when $F$ is Fano with $\mathscr{L}_F = \mathcal{O}_F(-K_F)$. \label{ex:SplittingFibreFano}
	\end{Example}
	In general, there might not exist an ample line bundle $\mathscr{L}_F$ which is invariant under $\Aut(F)$. The key in our situation is that we do not need to consider the full automorphism group:
	\begin{Lemma}
		Let $f\colon X\rightarrow Y$ be a proper holomorphic fibre bundle with fibre $F$ and denote $G:=\Aut(F)$. Assume that $X, Y$ and $F$ are normal varieties. If $f$ is projective, then the structure group of the fibre bundle $f$ can be reduced to an \emph{algebraic} subgroup $G'$ satisfying $G^0\subseteq G'\subseteq G$. Here, $G^0$ denotes the connected component of the identity of $G$.\label{reductionStructureGroupProjBundle}
	\end{Lemma}
	\begin{proof}
		Note to begin with that since $F$ is projective $G^0$ is certainly algebraic \cite[Theorem 2.3]{brion_AutomorphismGroups}. Then, any other group $G^0\subseteq G'\subseteq G$ is algebraic if and only if $G'/G^0$ is finite.
		
		As before, we denote by $p\colon \mathcal{G}= \mathcal{G}(f)\rightarrow Y$ the holomorphic principal $G$-bundle constructed using the same transition functions as $f$. Let us consider the complex analytic variety $Y':=\mathcal{G}/G^0\overset{\pi}{\rightarrow}Y$. Then, $\pi$ is a topological covering (with $Y'$ possibly non-connected) with transitively acting covering group $G/G^0$. In particular, $\pi$ is defined by a homomorphism 
		\begin{align}
		\rho\colon \pi_1(Y)\rightarrow G/G^0.\label{eqs:splittingTangentSequence3}
		\end{align}
		Let us denote the image of this map by $K$. Certainly the structure group of $f$ can be reduced to the group $G':= pr^{-1}(K)$, where $pr\colon G\rightarrow G/G^0$ is the canonical map. By construction $G^0\subseteq G'\subseteq G$ and it remains to prove that $K = G'/G^0$ is finite.
		
		To this end we consider the $\Z$-local system $\E:=\mathcal{G}\times^G \textmd{N}^1(F)$ on $Y$. Note that the natural action of $G=\Aut(F)$ on $\textmd{N}^1(F)$ factors through $G/G^0$. It follows that $\E$ is given by the representation $\pi_1(Y)\overset{\rho}{\rightarrow} G/G^0 \rightarrow \GL(\textmd{N}^1(F))$. Now, for any $f$-ample line bundle $\mathscr{L}$ on $X$ the class $[\mathscr{L}|_F]\in \textmd{N}^1(F)$ is $\pi_1(Y)$-invariant. In other words, the image $K\subseteq G/G^0$ of $\rho$ is contained in the stabiliser of the ample class $[\mathscr{L}|_F]$. But by \cite[Theorem 2.10]{brion_AutomorphismGroups} this stabiliser is finite. We deduce that $K$ is finite. By the argument following (\ref{eqs:splittingTangentSequence3}) we conclude that the structure group of $f$ can be reduced to the algebraic group $G':=pr^{-1}(K)$ as required.	
	\end{proof}
	Let us denote by $\mathcal{G}'$ the principal $G'$-bundle obtained by reduction of the structure group to $G'$ from $\mathcal{G}= \mathcal{G}(f)$.  
	\begin{Proposition}
		Let $f\colon X\rightarrow Y$ be a holomorphic fibre bundle with fibre $F$. Assume that $X, Y$ and $F$ are normal varieties. Assume, moreover, that $f$ is projective and that $\mathrm{Pic}^0(F) = 0$. Then, using the notation introduced above, there exists a holomorphic line bundle $\mathscr{L}$ on $X$ satisfying the following properties:
		\begin{itemize}
			\item[(1)] $\mathscr{L}$ is $f$-relatively ample and
			\item[(2)] $f_*\mathscr{L}^{\otimes m}$ is associated to $\mathcal{G}'$ for any $m>0$.
		\end{itemize}\label{existenceAssociatedAmpleBundle}
	\end{Proposition}
	\begin{proof}
		Let $\mathscr{L}'_F$ be any ample line bundle on $F$. Recall that $G'/G^0$ is finite by construction and pick a line bundle $\mathscr{L}_F$ in the numerical class of
		\begin{align*}
		\left[\mathscr{L}_F \right] = \sum_{g \in G'/G^0} g^*\left[\mathscr{L}'_F \right]\in N^1(F).
		\end{align*}
		Then $\mathscr{L}_F$ is ample as a linear combination of ample classes and clearly its numerical class is $G'$-invariant. Possibly replacing $\mathscr{L}_F$ by some multiple (to get rid of the torsion) we can assume that even
		\begin{align*}
		    \left[\mathscr{L}_F \right] \in \mathrm{Pic}(F) = \mathrm{Pic}(F)/\mathrm{Pic}^0(F) \hookrightarrow \mathrm{H}^2(F, \mathds{Z})
		\end{align*}
		is $G'$-invariant. By Lemma \ref{linearisationsAmpleBundles} some tensor power $\mathscr{L}_F^{\otimes k}$ admits a $G'$-linearisation. Replace $\mathscr{L}_F$ by this tensor power and let us denote $\mathscr{L} := \mathcal{G}'\times^{G'} \mathscr{L}_F$. Then $\mathscr{L}$ can be considered as a holomorphic line bundle on $X$, $\mathscr{L}|_F = \mathscr{L}_F$ is ample and
		\begin{align*}
		f_*\mathscr{L}^{\otimes m} = \mathcal{G}'\times^{G'} \mathrm{H}^0\Big(F, \mathscr{L}_F^{\otimes m}\Big) 
		\end{align*}
		 is associated to $\mathcal{G}'$ for any $m$.
	\end{proof}
	
	\subsection{Step 2: Existence of a Holomorphic Connection}
	
	Throughout this subsection we again fix a proper holomorphic fibre bundle $f\colon X\rightarrow Y$ and we assume that $Y$ is smooth, that the fibre $F$ is a normal variety and that the relative tangent bundle sequence 
	\begin{align}
	0\rightarrow \mathcal{T}_{X/Y} \rightarrow \mathcal{T}_X\rightarrow f^*\mathcal{T}_Y \rightarrow 0\label{eqs:splittingTangentSequence1}
	\end{align}
	splits. We want to show that the principal $G:=\Aut(F)$-bundle $\mathcal{G}:= \mathcal{G}(f)$ admits a holomorphic connection. In other words, we want to show that the sequence 
	\begin{align}
	0 \rightarrow \mathfrak{ad}(\mathcal{G}) \rightarrow \left(p_*\mathcal{T}_{\mathcal{G}}\right)^G \rightarrow \mathcal{T}_Y \rightarrow 0\label{eqs:splittingTangentSequence1.5}
	\end{align}
	admits a splitting, cf.\ Remark \ref{rem:FlatPrincipalBundles}. This will be done by relating the sequences (\ref{eqs:splittingTangentSequence1}) and (\ref{eqs:splittingTangentSequence1.5}). Pushing down (\ref{eqs:splittingTangentSequence1}) by $f$ we obtain the exact sequence
	\begin{align}
	0\rightarrow f_*\mathcal{T}_{X/Y} \rightarrow f_*\mathcal{T}_X\rightarrow f_*f^*\mathcal{T}_Y.
    \label{eqs:splittingTangentSequence1.75}
	\end{align}
	Note that $f_*f^*\mathcal{T}_Y \cong \mathcal{T}_Y\otimes f_*\mathcal{O}_X \cong \mathcal{T}_Y$ thanks to the projection formula and Zariski's main theorem. Additionally, (\ref{eqs:splittingTangentSequence1.75}) is exact on the right and even split as the push-forward of any splitting of (\ref{eqs:splittingTangentSequence1}) will be a splitting for (\ref{eqs:splittingTangentSequence1.75}).
	
	In summary, we find the split short exact sequence
	\begin{align}
	0\rightarrow f_*\mathcal{T}_{X/Y} \rightarrow f_*\mathcal{T}_X\rightarrow \mathcal{T}_Y \rightarrow 0.\label{eqs:splittingTangentSequence2}
	\end{align}
	\begin{Proposition}
		Let $f\colon X \rightarrow Y$ be a proper holomorphic fibre bundle with fibre $F$. Assume that $F$ is a normal variety and that $Y$ is a smooth variety. Then the two short exact sequences
		\begin{align*}
		0\rightarrow f_*\mathcal{T}_{X/Y} \rightarrow f_*\mathcal{T}_X\rightarrow \mathcal{T}_Y \rightarrow 0\\
		0 \rightarrow \mathfrak{ad}(\mathcal{G}) \rightarrow \left(p_*\mathcal{T}_{\mathcal{G}}\right)^G \rightarrow \mathcal{T}_Y \rightarrow 0
		\end{align*}
		can be naturally identified. In particular, if (\ref{eqs:splittingTangentSequence1}) splits then also (\ref{eqs:splittingTangentSequence1.5}) splits so that $\mathcal{G}$ admits a holomorphic connection.\label{existenceOfConnectionSplitTangentBundle}
	\end{Proposition}
	\begin{proof}
        Indeed, it suffices to show that the sequences can be identified locally in a natural way and that this identification is compatible with transition functions. To this end, let us first suppose that $X = U\times F\rightarrow U$ is the trivial fibration, where $U\subseteq \mathds{C}^n$ is some open subset. Then the sequences above can be naturally identified with
		\begin{align*}
		\begin{matrix}
		0 & \rightarrow & \mathrm{H}^0(F, \mathcal{T}_F)\otimes_{\mathds{C}} \mathcal{O}_U & \rightarrow & \mathrm{H}^0(F, \mathcal{T}_F)\otimes_{\mathds{C}} \mathcal{O}_U \oplus \mathcal{T}_U & \rightarrow & \mathcal{T}_U \rightarrow 0\\
		0 & \rightarrow & \mathfrak{g}\otimes_{\mathds{C}} \mathcal{O}_U & \rightarrow & \mathfrak{g}\otimes_{\mathds{C}} \mathcal{O}_U\oplus \mathcal{T}_U & \rightarrow & \mathcal{T}_U \rightarrow 0.
		\end{matrix}
		\end{align*}
		Since there exist natural identifications $\mathrm{H}^0(F, \mathcal{T}_F) = \mathfrak{g} = \mathfrak{g}(\Aut(F))$, see e.g.\ \cite[Remark 2.4.(iii)]{brion_AutomorphismGroups}, this proves the assertion in case $X = U\times F$. 
		
		Now, let $\phi\colon U\rightarrow G = \Aut(F)$ be a holomorphic function. We need to verify that the above identifications are compatible with the isomorphism of fibre bundles $U\times F \rightarrow U\times F, (x, z) \mapsto (x, \phi_x(z))$ defined by $\phi$. But indeed, for any fixed $x\in U$ the induced transition functions for $f_*\mathcal{T}_X$, respectively $\left(p_*\mathcal{T}_{\mathcal{G}}\right)^G$, are:
		\begin{align*}
		\mathrm{H}^0(F, \mathcal{T}_F) \oplus \mathcal{T}_U|_x \rightarrow \mathrm{H}^0(F, \mathcal{T}_F) \oplus \mathcal{T}_U|_x, \quad 
		(V, w) \mapsto \left( \begin{matrix}
		d\phi_x\circ \phi^{-1} & d_x\phi \\ 0 & \mathrm{id}
		\end{matrix} \right) 
		\left( \begin{matrix}
		V \\ w
		\end{matrix} \right) \\
		\mathfrak{g} \oplus \mathcal{T}_U|_x \rightarrow \mathfrak{g} \oplus \mathcal{T}_U|_x, \qquad \qquad 
		(v, w) \mapsto \left( \begin{matrix}
		\mathrm{Ad}_{\phi_x} & d_x\ell_{\phi}|_1 \\ 0 & \mathrm{id}
		\end{matrix} \right) 
		\left( \begin{matrix}
		v \\ w
		\end{matrix} \right)\hspace{0.3cm}
		\end{align*}
		Using again the identification $\mathrm{H}^0(F, \mathcal{T}_F) = \mathfrak{g} = \mathfrak{g}(\Aut(F))$ we see that these are the same and so we conclude.
	\end{proof}

	\subsection{Step 3: Conclusion of the Proof}
	
	We are now ready to finish the proof of Lemma \ref{splittingTangentSequence}: Consider the principal $\Aut(F)$-bundle $\mathcal{G}(f)$. From Proposition \ref{existenceOfConnectionSplitTangentBundle} $\mathcal{G}|_{Y_{\mathrm{reg}}}$ inherits a holomorphic connection. By Lemma \ref{reductionStructureGroupProjBundle} there exists an algebraic subgroup $\Aut^0(F)\subseteq G'\subseteq \Aut(F)$ so that $\mathcal{G}(f)$ admits a reduction to a $G'$-bundle $\mathcal{G}'$. Note that also $\mathcal{G}'|_{Y_{\mathrm{reg}}}$ admits a holomorphic connection: Indeed, there are two ways to see this. Either, one notes that the proof of Lemma \ref{reductionStructureGroupProjBundle} shows that the reduction from $\mathcal{G}$ to $\mathcal{G}'$ is purely topological and, hence, is compatible with any connection. Alternatively, one notes that the sequence (\ref{eqs:splittingTangentSequence1.5}) clearly agrees with the corresponding one for $G'$ as $\mathfrak{g}=\mathfrak{g}'$. 
	
	Now, fix an $f$-ample line bundle $\mathscr{L}$ on $X$ as in Proposition \ref{existenceAssociatedAmpleBundle} so that $f_*\mathscr{L}^{\otimes m}$ is associated to $\mathcal{G}'$ for any $m>0$. Since $\mathcal{G}'|_{Y_{\mathrm{reg}}}$ admits a holomorphic connection the same is true of its associated bundle $f_*\mathscr{L}^{\otimes m}|_{Y_{\mathrm{reg}}}$. By Proposition \ref{flatnessOnCYs} we conclude that $f_*\mathscr{L}^{\otimes m}$ is numerically flat for any $m>0$. Thus, $f$ is locally constant by the by-now standard criterion \cite[Proposition 2.1]{wang_manifoldsWithNefAnticanonicalBundle} (to the best of the author's knowledge the latter argument first appeared in \cite[Proposition 4.3.6]{cao_PhDThesis}).

	\section{Locally Constant Fibrations and Positivity of Curvature}
	\label{sec:LocConstFibrAndPositiveCurvature}
	
	In this section we study when the total space of a locally constant fibre bundle is positively curved. This will complete, in particular, the proof of Theorem \ref{IntroCharacterisationMfdsNefACBundle}.

    \subsection{Varieties with Nef Anti-Canonical Bundle}
    
	\begin{Theorem}
		Let $Y$ be a log terminal projective variety with numerically trivial canonical class $\cc_1(Y) = 0$, let $(F, \Delta_F)$ be a projective klt pair such that $-(K_F + \Delta_F)$ is nef and let $\rho\colon \pi_1(Y) \rightarrow \Aut(F, \Delta_F)$ be any group homomorphism. Let us denote by $\widetilde{Y}\rightarrow Y$ the universal cover of $Y$ and set
		\begin{align*}
		    (X, \Delta) := \left(\widetilde{Y}\times F,\ pr_F^*\Delta_F \right)/\pi_1(Y).
		\end{align*}
        Then the pair $(X, \Delta)$ has klt singularities. Moreover, if $X$ is projective then the anti-log canonical divisor $-(K_X+ \Delta)$ is nef.\label{charMfdsNefALogCanBundle}
	\end{Theorem}
    Conversely, it follows from the successive works \cite{paun_FGNefACBundle}, \cite{paun_FGNefACBundleII}, \cite{zhang_VarietiesNefACBundle}, \cite{lu_SemiStabilityOfAlbanese}, \cite{cao_PhDThesis},  \cite{cao_NefAnticanonicalBundle}, \cite{cao_NefAnticanonicalBundleII}, \cite{cao_NefAnticanonicalBundleIII}, \cite{CCM_VarietiesNefACBundle}, \cite{wang_manifoldsWithNefAnticanonicalBundle} and \cite{matsumuraWang_NefAnticanonicalBundle} that, up to finite quasi-\'etale cover, any projective klt pair $(X, \Delta)$ such that $-(K_X+ \Delta)$ is nef admits a locally constant fibration $f\colon X \rightarrow Y$ as in Theorem \ref{charMfdsNefALogCanBundle}.
	\begin{proof}
		First of all, $(X, \Delta)$ has klt singularities since this property is local in the analytic topology and since the pair $(\widetilde{Y}\times F,\ pr_F^*\Delta_F)$ is klt.
  
  Regarding the second assertion, let us assume that $X$ is projective and denote by $f\colon X\rightarrow Y$ the natural morphism. By Lemma \ref{criterionForProjectivity}, there exists a $\pi_1(Y)$-invariant, linearisable very ample line bundle $\mathscr{L}_F$ on $F$. Then, clearly, $\mathscr{L} := (\widetilde{Y}\times \mathscr{L}_F)/\pi_1(Y)$ is an $f$-relatively very ample holomorphic line bundle on $X$.
		
		Now, for any fixed integer $m>0$ and any $\ell>0$ we consider the coherent $\mathcal{O}_Y$-modules
		\begin{align*}
		 f_*&\left(\left(\mathscr{L}\otimes\mathcal{O}_X\left(-m\left(K_{X/Y}+\Delta\right)\right)\right)^{\otimes \ell}\right) \\
		 &\qquad = \left( \widetilde{Y} \times \mathrm{H}^0\Big(F, \big( \mathscr{L}_F\big(-m \left(K_{F}+\Delta_F\right)\big) \big)^{\otimes \ell} \Big) \right) /\pi_1(Y).
		\end{align*}
		Observe that the $f_*((\mathscr{L}\otimes\mathcal{O}_X(-m(K_{X/Y}+\Delta)))^{\otimes \ell})$ are holomorphically flat and, hence, by Proposition \ref{flatnessOnCYs} also numerically flat vector bundles on $Y$. The rest of the argument is standard following \cite[Proof of Theorem 4.7]{matsumura_asymptoticBaseLoci}: Since numerically flat bundles are nef and, since nefness is preserved by pull-back \cite[Proposition 6.1.8]{lazarsfeld_PositivityII}, also
		\begin{align}
		    f^*f_*\left(\left(\mathscr{L}\otimes\mathcal{O}_X\left(-m\left(K_{X/Y}+\Delta\right)\right)\right)^{\otimes \ell}\right).\label{eqs:CharMfdNefACBundle_1}
		\end{align}
		is nef. Note that the line bundle $\mathscr{L}_F\otimes\mathcal{O}_F(-m(K_{F}+\Delta_F))$ is ample for any $m>0$ as the product of an ample bundle and a nef bundle. Consequently,
		\begin{align*}
		\Big(\mathscr{L}_F\otimes\mathcal{O}_F\left(-m\left(K_{F}+\Delta_F\right)\right)\Big)^{\otimes \ell}
		\end{align*}
		is generated by global sections for any sufficiently large $\ell = \ell(m)\gg0$. It follows, that the natural morphism
		\begin{align*}
		f^*f_*\left(\left(\mathscr{L}\otimes\mathcal{O}_X\left(-m\left(K_{X/Y}+\Delta\right)\right)\right)^{\otimes \ell}\right) \twoheadrightarrow 
		\left( \mathscr{L}\otimes\mathcal{O}_X\left(-m\left(K_{X/Y}+\Delta\right)\right) \right)^{\otimes \ell}
		\end{align*}
		is surjective for all sufficiently large $\ell\gg0$. Since the bundles on the left are nef according to the discussion above (\ref{eqs:CharMfdNefACBundle_1}) and, since quotients of nef bundles are nef \cite[Theorem 6.2.12]{lazarsfeld_PositivityII}, we conclude that also $\big(\mathscr{L}\otimes\mathcal{O}_X(-m(K_{X/Y}+\Delta))\big)^{\otimes \ell}$ and, hence, $\mathscr{L}\otimes\mathcal{O}_X(-m(K_{X/Y}+\Delta))$ are nef for any $m>0$. Since limits of nef bundles are nef \cite[Proposition 6.2.11]{lazarsfeld_PositivityII}, we deduce that also $-(K_{X/Y}+\Delta)$ is nef. But then so is
		\begin{align*}
			-(K_X + \Delta) = -(K_{X/Y} + \Delta + f^*K_Y) \equiv -(K_{X/Y} + \Delta).
		\end{align*}
		Here, we use that $K_Y \equiv 0$ is numerically trivial.
	\end{proof}
    \begin{Remark}
        Let $\mathscr{L}$ be as in the proof of Theorem \ref{charMfdsNefALogCanBundle}. Then the line bundle 
        $$
        \mathscr{L}_{m, \ell} := \left(\mathscr{L}\otimes\mathcal{O}_X\left(-m\left(K_{X/Y}+\Delta\right)\right)\right)^{\otimes \ell}
        $$
        is nef and $f$-relatively very ample for any $m\geq1$ and any $\ell = \ell(m) \gg 1$. Moreover, $f_*\mathscr{L}_{m, \ell}$ is numerically flat and associated to $\rho$.\label{rem:relatively ample bundle}
    \end{Remark}
	\begin{Example}
	    Let $Y$ be a log terminal projective variety with numerically trivial canonical class $\cc_1(Y)=0$, let $F$ be a (weak) Fano and let $\rho\colon \pi_1(Y)\rightarrow \Aut(F)$ be any group homomorphism. As before, we denote by $X:= (\widetilde{Y}\times F)/\pi_1(Y)$ the total space of the associated locally constant fibration $f\colon X\rightarrow Y$. Then by Corollary \ref{criterionForProjectivitySpecialCase} and Theorem \ref{charMfdsNefALogCanBundle} $X$ is projective and the anti-canonical divisor $-K_X$ is nef.
	\end{Example}
	In particular, we see how Theorem \ref{charMfdsNefALogCanBundle} allows to explicitly construct many new and non-trivial examples of varieties with nef anti-log canonical bundle.
    \begin{proof}[Proof of Theorem \ref{IntroCharacterisationMfdsNefACBundle}]
        The assertion on the projectivity is just Lemma \ref{criterionForProjectivity} while the assertion that the anti-canonical bundle is nef follows from Theorem \ref{charMfdsNefALogCanBundle}.
    \end{proof}
	One might wonder what happens in the situation of Theorem \ref{IntroCharacterisationMfdsNefACBundle} in case one asks for more positivity of the anti-canonical bundle. The forward implications in the following statement are (with minor modifications) due to Campana-Demailly-Peternell \cite{campanaDemaillyPeternell_SemiPositiveRicci} and Ejiri-Iwai-Matsumura \cite{matsumura_asymptoticBaseLoci} and they are clearly sharp.
    \begin{Proposition}\emph{(c.f.\ \cite{campanaDemaillyPeternell_SemiPositiveRicci, matsumura_asymptoticBaseLoci})}

        \noindent
        Let $Y$ be a log terminal projective variety with numerically trivial canonical class $\cc_1(Y) = 0$, let $(F, \Delta_F)$ be a projective klt pair such that $-(K_F + \Delta_F)$ is nef and let $\rho\colon \pi_1(Y) \rightarrow \Aut(F, \Delta_F)$ be any group homomorphism. Let us denote by $\widetilde{Y}\rightarrow Y$ the universal cover of $Y$ and set
		\begin{align*}
		    (X, \Delta) := \left(\left(\widetilde{Y}\times F\right),\ pr_F^*\Delta_F \right)/\pi_1(Y).
		\end{align*}
        Assume that $X$ is projective. Then
        \begin{itemize}
            \item[(1)] If $-(K_X + \Delta)$ is Hermitian semipositive then the image of $\rho$ is contained in a compact subgroup for the Euclidean topology. 
            \item[(2)] Conversely, if the image of $\rho$ is contained in a compact subgroup for the Euclidean topology  and if $-(K_F + \Delta_F)$ is Hermitian semipositive then also $-(K_X + \Delta)$ is Hermitian semipositive.
            \item[(3)] If $-(K_X + \Delta)$ is semiample then the image of $\rho$ is finite. 
            \item[(4)] Conversely, if the image of $\rho$ is finite and if $-(K_F + \Delta_F)$ is semiample then also $-(K_X + \Delta)$ is semiample.
        \end{itemize}\label{charOtherPosNotions}
    \end{Proposition}
    \begin{proof}
        Let $\mathscr{L}_{m, \ell}$ be defined as in Remark \ref{rem:relatively ample bundle}. Let $\textmd{Im}(\rho) \subseteq G \subseteq \Aut(F)$ denote the Zariski-closure of $\textmd{Im}(\rho)$. By Lemma \ref{criterionForProjectivity}, $G$ is linear algebraic. After possibly replacing $X$ by some finite \'etale cover we can also assume that $G$ is connected.
        
        If $-(K_X + \Delta)$ admits a smooth Hermitian metric of semipositive curvature then $f_*\mathscr{L}_{m, \ell}$
        admits a singular Hermitian metric of semipositive curvature \cite{BP_PositivityDirectImages}, see \cite[Proof of Theorem 4.8]{matsumura_asymptoticBaseLoci} for the precise argument. By \cite[Proposition 3.3]{IwaiMatsumuraZhong_SingularPositiveTangentSheaf}, the vector bundle $f_*\mathscr{L}_{m, \ell}$ is unitarily flat, i.e. the image of
        $$
        \pi_1(Y) \overset{\rho}{\rightarrow} \Aut^0(F) \rightarrow \GL\left( \mathrm{H}^0\Big(F, \mathscr{L}_{m, \ell}|_F \Big) \right)
        $$
        is contained in a Euclidian compact subgroup. Since the natural group homomorphism $G \rightarrow \GL(\mathrm{H}^0(F, \mathscr{L}_{m, \ell}|_F))$ is injective (c.f.\ the proof of Lemma \ref{linearisationsAmpleBundles}) it follows that also $\textmd{Im}(\rho) \subseteq \Aut(F)$ is contained in a Euclidian compact subgroup. Conversely, if the image of $\rho$ is contained in a Euclidean compact subgroup $K \subseteq G$ then, by averaging over the $K$-action, we can find a $K$-invariant smooth Hermitian metric of semipositive curvature on $-(K_F + \Delta_F)$. Clearly, this metric extends to a metric on $-(K_{X/Y} + \Delta)$ with the same property.
        
        Let us now assume that $-(K_X + \Delta)$ is semiample. Following \cite[Proof of Theorem 1.7]{matsumura_asymptoticBaseLoci} we consider the anti-Iitaka fibration
        $$
        \phi\colon F \rightarrow P_F := \textmd{Proj}\left( \bigoplus_{m=0}^\infty \mathrm{H}^0\Big(F, \big( \mathcal{O}_F\big(-m \left(K_{F}+\Delta_F\right)\big)\big)\Big) \right).
        $$
        Then, as in \cite{matsumura_asymptoticBaseLoci}, the action of $G$ on $P_F$ is trivial. The reason is that the induced morphism $p\colon P_F \rightarrow P_X = \textmd{Proj}(\oplus \mathrm{H}^0(X, ( \mathcal{O}(-m (K_{X}+\Delta)))))$ is finite \cite[Proof of Theorem 1.7, Step 2]{matsumura_asymptoticBaseLoci} and that the action of $\pi_1(Y)$ on $P_X$ is trivial \cite[Proof of Theorem 1.7, Step 3]{matsumura_asymptoticBaseLoci}. If follows that $G$ acts effectively on a general fibre $F_p$ of $\phi$. By construction, $(F_p, \Delta|_{F_p})$ is a klt pair and $K_{F_p} + \Delta|_{F_p} \sim_{\mathds{Q}} 0$. Then, by \cite[Theorem 4.5]{Xu_FibrationsOnLogCYs}, $\Aut^0(F_p, \Delta|_{F_p})$ is an abelian variety. We conclude that, as a connected linear algebraic subgroup of an abelian variety, $G = 1$ is trivial, i.e.\ that $(X, \Delta) \cong Y \times (F, \Delta_F)$. The converse implication is obvious.
    \end{proof}

    Note that the reverse implications $(2), (4)$ in Proposition \ref{charOtherPosNotions} are immediate and do not make use of our assumptions on $Y$. In particular, the corresponding statements for the relative anti-canonical divisor $-(K_{X/Y} + \Delta)$ remain valid over any (normal, projective) base $Y$. On the other hand, the conclusion that $-K_X \equiv -K_{X/Y}$ is nef in Theorem \ref{IntroCharacterisationMfdsNefACBundle} really is non-trivial as the following example shows:
    \begin{Example}
        Let $C$ be a smooth projective curve of genus $2$ and write $K_C = \mathcal{O}_C(p + p')$. A quick computation shows that there exists a unique non-split extension
        $$
        0 \rightarrow \mathcal{O}_C(p) \rightarrow \mathcal{E} \rightarrow \mathcal{O}_C(-p') \rightarrow 0.
        $$
        Then $\mathcal{E}$ is a holomorphic vector bundle of $\rk\mathcal{E} = 2$ and $\deg(\mathcal{E}) = 0$. Moreover, it is clearly not nef. We claim that $\mathcal{E}$ can not be written $\mathcal{E} = \mathscr{L}_1 \oplus \mathscr{L}_2$ as the direct sum of two line bundles. Indeed, if this were true then the two maps
        $$
        \mathscr{L}_i \hookrightarrow \mathcal{E} \rightarrow \mathcal{O}_C(-p')
        $$
        would be non-trivial, hence $\deg(\mathscr{L}_i) \leq -1$ for $i=1,2$, a contradiction. Thus, by \cite[Proposition 19]{atiyah_ConnectionsInFibreBundles}, $\mathcal{E}$ is holomorphically flat, given by a representation $\rho\colon \pi_1(C) \rightarrow \GL_2(\C)$ say.

        Now, consider $X := (\mP^1 \times C)/\pi_1(C) = \mP(\mathcal{E})$. Then
        $$
        \mathcal{O}_X(-K_{X/C}) 
        = \mathcal{O}_{\mP(\mathcal{E})}(2) 
        $$
        is not nef. Observe that this is so, precisely because of the failure of Proposition \ref{IntroBiswas} for $Y = C$ and $\mathcal{E}$ as above.
        \label{ex:PositivityCanFail}
    \end{Example}

    \subsection{Varieties with Nef Tangent Bundle}

	\begin{Proposition}
		Let $T$ be an abelian variety of dimension $q$, let $F$ be a smooth Fano variety with nef tangent bundle and let $\rho\colon \pi_1(T) \cong \Z^{2q} \rightarrow \Aut(F)$ be a group homomorphism. Set $X:= (F \times \C^q)/\pi_1(T)$. Then $X$ is a smooth projective variety and the tangent bundle of $X$ is nef.\label{characterisationMfdNefTangent}
	\end{Proposition}
	Conversely, it was proved by Demailly-Peternell-Schneider \cite[Main Theorem]{DPS_ManifoldsWithNefTangentBundle} and Cao \cite{cao_PhDThesis} that, up to finite \'etale cover, every smooth projective variety with nef tangent bundle is of the form in Proposition \ref{characterisationMfdNefTangent}.
	\begin{proof}
		That $X$ is projective follows from Corollary \ref{criterionForProjectivitySpecialCase} as $F$ is Fano. The rest of the argument is essentially identical to the one in Theorem \ref{charMfdsNefALogCanBundle}.
	\end{proof}

	\bibliographystyle{alpha}
\bibliography{LocConst.bib}
\end{document}